\documentclass[leqno,12pt]{article}

\usepackage{amsfonts}

\evensidemargin\oddsidemargin

\newcommand{\eqnsection}{
\renewcommand{\theequation}{\thesection.\arabic{equation}}
    \makeatletter
    \csname  @addtoreset\endcsname{equation}{section}
    \makeatother}
\eqnsection

\newcommand{\e}{{ \rm  {e}}}

\renewcommand{\theequation}{\thesection.\arabic{equation}}
\newtheorem{theorem}{Theorem}[section]
\newtheorem{lemma}[theorem]{Lemma}
\newtheorem{proposition}[theorem]{Proposition}
\newtheorem{corollary}[theorem]{Corollary}

\begin{document}

\vglue15pt
\centerline{\LARGE Tail asymptotics for the total progeny of the }

\bigskip
\centerline{\LARGE critical killed branching random walk}

\bigskip\bigskip
\centerline{Elie A\"id\'ekon\footnote{ Eurandom, Technische Universiteit Eindhoven, P.O. Box 513, 5600 MB Eindhoven, The Netherlands. email: elie.aidekon@gmail.com}}
\bigskip

\bigskip
\bigskip

{\leftskip=2truecm \rightskip=2truecm \baselineskip=15pt \small

\noindent{\slshape\bfseries Summary.} We consider a branching random walk on $\mathbb{R}$ with a killing barrier at zero. At criticality, the process becomes eventually extinct, and the total progeny $Z$ is therefore finite. We show that  $P(Z>n)$ is of order $(n\ln^2(n))^{-1}$, which confirms the prediction of Addario-Berry and Broutin \cite{addario}.



}

\section{Introduction}

We look at the branching random walk on $\mathbb{R}_+$ killed below zero. Let $b\ge 2$ be a determinist integer which represents the  number of children of the branching random walk, and $x\ge 0$ be the position of the (unique) ancestor. We introduce the rooted $b$-ary tree $\mathcal{T}$, and we attach at every vertex $u$ except the root an independent random variable $X_u$ picked from a common distribution (we denote by $X$ a generic random variable having this distribution). We define the position of the vertex $u$ by
	$$
	S(u) := x + \sum_{v<u}X_v
	$$
where $v<u$ means that the vertex $v$ is an ancestor of $u$. We say that a vertex (or particle) $u$ is alive if $S(v)\ge 0$ for any ancestor $v$ of $u$ including itself. \\

The process can be seen in the following way. At every time $n$, the living particles split into $b$ children. These children make independent and identically distributed steps. The children which enter the negative half-line are immediately killed and have no descendance.  We are interested in the behaviour of the surviving population. At criticality (see below for the definition), the population ultimately dies out. We define the total progeny $Z$ of the killed branching random walk by
	$$
	Z:=\#\{ u\in \mathcal{T}\,:\, S(v)\ge 0 \; \forall \, v\le u\}\,.
	$$
Aldous \cite{aldous} conjectured that in the critical case, $E[Z]<\infty$ and $E[Z\ln(Z)]=\infty$. In \cite{addario}, Addario-Berry and Broutin proved that conjecture (in a more general setting where the number of children may be random). As stated there, this is a strong hint that $P(Z=n)$ behaves asymptotically like $1/(n^2\ln^2(n))$, which is a typical behaviour of critical killed branching random walks. Here, we look at the tail distribution $P(Z\ge n)$. We mention that the Branching Brownian Motion, which can be seen as a continuous analogue of our model, already drew some interest. Kesten \cite{kest78} and Harris and Harris \cite{harris} studied the extinction time of the population, whereas Berestycki et al. \cite{Berestgenealogy} showed a scaling limit of the process near criticality. Maillard \cite{maillard} investigated the tail distribution of $Z$, and proved that $P(Z=n) \sim {c \over n^2 \ln^2 n}$ as expected. \\

Before stating our result, we introduce the Laplace transform $\phi(t) := E[\e^{tX}]$ and we suppose that
\begin{itemize}
\item $\phi(t)$ reaches its infimum at a point $t=\rho>0$ which belongs to the interior of $\{t\,:\phi(t)<\infty\}$,
\item The distribution of $X$ is non-lattice.
\end{itemize}

\noindent The second assumption is for convenience in the proof, but the theorem remains true in the lattice case. The probability that the population lives forever is zero or positive depending on whether $E[e^{\rho X}]$ is less or greater than the critical value $1/b$. In the present work, we consider the critical branching random walk which corresponds to the case $E[e^{\rho X}]=1/b$. For $x\ge 0$, we call $P^x$ the distribution of the killed branching random walk starting from $x$.

\begin{theorem}\label{main}
There exist two positive constants $C_1$ and $C_2$ such that for any $x\ge 0$, we have for $n$ large enough
	$$
	C_1{(1+x)\e^{\rho x}\over n\ln^2(n)}\le P^x(Z>n) \le C_2 {(1+x)\e^{\rho x}\over n\ln^2(n)}\,.
	$$
\end{theorem}

\bigskip

Hence, the tail distribution has the expected order. Nevertheless, the question to find an equivalent to $P(Z=n)$ is still open. As observed in \cite{addario}, in order to have a big population, a particle of the branching random walk needs to go far to the right, so that its descendance will be greater than $n$ with probability large enough (roughly a positive constant). The theorem then comes from the study of the tail distribution of the maximum of the killed branching random walk. By looking at the branching random walk with two killing barriers, we are able to improve the estimates already given in \cite{addario}.

\bigskip

The paper is organised as follows. Section 2 gives some elementary results for one-dimensional random walks on an interval. Section 3 gives estimates on the first and second moments of the killed branching random walk, while Section 4 contains the asymptotics on the tail distribution of the maximal position reached by the branching random walk before its extinction. Finally, Theorem \ref{main} is proved in Section 5.

\section{Results for one-dimensional random walks}

Let $ R_n = R_0 + Y_1 + \ldots + Y_n $ be a one-dimensional random walk and $P^x$ be the distribution of the random walk starting from $x$. For any $k\in \mathbb{R}$, we define $\tau_k^+$ (resp. $\tau_k^-$) as the first time the walk hits the domain $(k,+\infty)$ (resp. $(-\infty,k)$),
    \begin{eqnarray*}
    \tau_k^+ &:=& \inf\{ n \ge 0\,: \, R_n > k \}\,, \\
    \tau_k^- &:=& \inf\{ n \ge 0\,: \, R_n < k \}\,.
    \end{eqnarray*}

We assume
$$
(H) 	\qquad
		E[Y_1] = 0, \; \exists \,\theta,\,\eta >0 \mbox{ such that }
		E[\e^{-(\theta+\eta) Y_1}]<\infty,\, E[\e^{(1+\eta)Y_1}]<\infty.
$$

All the results of this section are stated under condition $(H)$. The results remain naturally true after renormalization as long as $E[e^{tY_1}]$ is finite on a neighborhood of zero (and $E[Y_1]=0$). Throughout the paper, the variables $C_1$, $C_2$, $\ldots$ represent positive constants. We first look at the moments of the overshoot $U_k$ and undershoot $L_k$ defined respectively by
	\begin{eqnarray*}
	U_k &:=& S_{\tau_k^+} - k \,,\\
	L_k &:=& k-S_{\tau_k^-}\,.
	\end{eqnarray*}

\begin{lemma}
\label{overshoot}
There exists $C_3>0$ such that $E^0[\e^{U_k}] \in [C_3,1/C_3]$ for any $k\ge 0$ and $E^0[\e^{\theta L_k}] \in [C_3,1/C_3]$ for any $k\le 0$.
\end{lemma}
{\it Proof}. This is a consequence of Proposition 4.2 in Chang \cite{chang}. $\Box$

\bigskip
The following lemma concerns the well-known hitting probabilities of $R$. 
\begin{lemma}
    \label{passage}
     For any $x\ge 0$,
        \begin{equation} \label{passage2}
        P^x(\tau_k^+<\tau_0^-) = { E[-S_{\tau_x^-}] \over k} + o(1/k)\,.
        \end{equation}
    as $k\to \infty$. Moreover, there exist two positive constants $C_4$ and $C_5$ such that, for any real $k\ge 0$ and any $z \in [0,k]$, we have
        \begin{equation}\label{passage1}
        C_4{z+ 1 \over k+1} \le P^z( \tau_k^+ <\tau_0^- ) \le C_5{z+1 \over k+1}\,.
        \end{equation}
\end{lemma}
{\it Proof}. Let $k> 0$ and $x \in[0,k]$. By Lemma \ref{overshoot}, we are allowed to use the stopping time theorem on $(R_n,\,n\le \min(\tau_0^-,\tau_k^+))$, and we get
    $$
    x = E^x[R_{\tau_k^+},\,\tau_k^+ < \tau_0^-] + E^x[R_{\tau_0^-},\, \tau_0^- < \tau_k^+]\,.
    $$
We can write it
    $$
    x = kP^x(\tau_k^+<\tau_0^-) + A_1 - A_2
    $$
where $A_1$ and $A_2$ are nonnegative and defined by $A_1 := E^x[U_k,\, \tau_k^+ < \tau_0^-]$ and $A_2 := E^x[L_0 ,\,\tau_0^- < \tau_k^+]$. Equivalently,
    \begin{equation}
    \label{passageeq}
    P^x(\tau_k^+ < \tau_0^-) = {x- A_1 +A_2 \over k}\,.
    \end{equation}
By Cauchy-Schwartz inequality and Lemma \ref{overshoot}, we observe that
    $$
    (A_1)^2 \le E^x[U_k^2]P^x(\tau_k^+<\tau_0^-) \le C_6 P^x(\tau_k^+<\tau_0^-)\,.
    $$
Since $P^x(\tau_k^+ < \tau_0^-)$ goes to zero when $k$ tends to infinity, we deduce that
    $$
    \lim_{k \to \infty} A_1 =0\,.
    $$
By dominated convergence, we have also
    $$
    \lim_{k \to \infty} A_2 = E^x[L_0]
    $$
and $E^x[L_0]\le C_7$ by Lemma \ref{overshoot}. This leads to equation (\ref{passage2}) since $E[-S_{\tau_x^-}] = x+E^x[L_0]$. Furthermore, we have
   $ 0 \le A_1 \le \sqrt{C_6}$ and
    $0\le A_2 \le C_7$. Therefore (\ref{passageeq}) implies that
    $$
    P^x(\tau_k^+ < \tau_0^-) \le {x + C_7 \over k} \le C_8{x+1 \over k+1}.
    $$
Similarly,
    $$
    P^x(\tau_k^+ < \tau_0^-) \ge {x-\sqrt{C_6} \over k}\,.
    $$
We notice also that $P^x(\tau_k^+ < \tau_0^-) \ge P^0(\tau_k^+ < \tau_0^-)$. By (\ref{passage2}), there exists a constant $C_9>0$ such that $P^0(\tau_k^+ < \tau_0^-) \ge {C_9 \over k+1}$. We get
    $$
    P^x(\tau_k^+ < \tau_0^-) \ge
    \left\{
        \begin{array}{cc}
        {C_9  \over k+1}  &\mbox{if } x< \sqrt{C_6} +1\\
        C_{10}{x+1\over k} &\mbox{otherwise}
        \end{array}
    \right.
    $$
with $C_{10} := {1\over \sqrt{C_6} +2}$. It implies that
$$
P^x(\tau_k^+ < \tau_0^-) \ge C_{11} {x+1\over k+1}.
$$
Thus equation (\ref{passage1}) holds with $C_5:= C_8$ and $C_4 := C_{11}$. $\Box$ \\

Throughout the paper, we will write $\Delta_k(1)$ for any function such that
	$$
	0 <D_1 \le \Delta_k(1) \le D_2
	$$

\noindent for some constants $D_1$ and $D_2$ and $k$ large enough. The following lemma provides us with estimates used to compute the moments of the branching random walk in Sections 3 and 4.

\begin{lemma}
\label{estimatepath}
We have for any $x>0$,
    \begin{eqnarray}
    \label{path0tok}
    E^0\left[ \e^{U_k}\sum_{\ell=0}^{\tau_{k}^+} \e^{-  R_{\ell}} (R_{\ell} + 1),\, \tau_{k}^+ <\tau_0^- \right] &=& \Delta_k(1){ 1 \over k}\,,\\
    \label{pathktok}
    E^{k-x}\left[\e^{U_k} \sum_{\ell=0}^{\tau_{k}^+} \e^{-  R_{\ell}} (R_{\ell} + 1),\, \tau_{k}^+ <\tau_0^- \right] &=& \Delta_k(1){ 1 +x \over k^2}\,,\\
    \label{pathkto0}
    E^{k-x}\left[\e^{-L_0} \sum_{\ell=0}^{\tau_{0}^-} \e^{-  R_{\ell}} (k-R_{\ell} + 1),\, \tau_{0}^- <\tau_k^+ \right] &=& \Delta_k(1)(1+x)\,.
    \end{eqnarray}
\end{lemma}
{\it Proof}. First let us explain how we can find intuitively these estimates. The terms of the sum within the expectation is big when $R_{\ell}$ is close to $0$, and the time that the random walk spends in the neighborhood of $0$ before hitting level $0$ is roughly a constant. Moreover, by Lemma \ref{overshoot}, we know that the overshoot $U_k$ and the undershoot $L_0$ behave like a constant. From here, we can deduce the different estimates. In (\ref{path0tok}), the optimal path makes the particle stay a constant time near zero then hit level $k$ which is of cost ${1/k}$. In (\ref{pathktok}), the particle first goes close to $0$, which gives a term in $(1+x)/k$, then go back to level $k$ which gives a term in $1/k$. Finally looking at (\ref{pathkto0}), we see that the particle goes directly to $0$, which brings a term of order $k$ because of the sum, and a term of order $(1+x)/k$ because of the cost to hit $0$ before $k$. The proofs of the three equations being rather similar, we restrain our attention on the proof of (\ref{path0tok}) for sake of concision.\\

\noindent We introduce the function $g(z):= \e^{- z} (1+z)$ and we observe that $g$ is decreasing. Let also
    $$
    A:=
    E^0\left[ \e^{U_k}\sum_{\ell=0}^{\tau_{k}^+} \e^{-  R_{\ell}} (R_{\ell} + 1),\, \tau_{k}^+ <\tau_0^- \right]\,.
    $$
Let $a>0$ be such that $P(Y_1>a)>0$ and $P(Y_1<-a)>0$. For ease of notation we suppose that we can take $a=1$. For any integer $i$ such that $0\le i<k$, we denote by $I_i$ the interval $[i,i+1)$, and we define
\begin{eqnarray*}
T_i  &:=&  \inf\{ n\ge 1\,:\, R_n \in I_i \}\,, \\
N(i) &:=&   \#\{ n\le \min\{\tau_k^+,\tau_0^-\}\,:\, R_n \in I_i \}
\end{eqnarray*}
which respectively stand for the first time the walk enters $I_i$ and the number of visits to the interval before hitting level $k$ or level $0$. We observe that
    \begin{eqnarray*}
    A  \le  \sum_{0\le i <k} g(i+1)E^0[\e^{U_k}N(i),\, \tau_{k}^+ <\tau_0^- ]\,.
    \end{eqnarray*}
Let $i$ be an integer between $1$ and $k-1$, and let $z \in [i,i+1)$. We have
	\[
	P^z(T_i>\min(\tau_0^-,\tau_k^+) ) \ge P^z(R_{\ell}\le R_1, \,\forall\, \ell\in [1,\tau_0^-],\,R_1< i)\,.
	\]
We use the Markov property to get
    \begin{eqnarray*}
    P^z( T_i>\min(\tau_0^-,\tau_k^+))
     \ge  E^z\left[P^{h}(\tau_{h}^+ >\tau_0^-)_{h=R_1},\,R_1<i\right] \,.
    \end{eqnarray*}
By Lemma \ref{passage} equation (\ref{passage2}) (applied to $-R$), there exists a positive constant $C_{12}$ such that 	 $P^{h}(\tau_{h}^+ >\tau_0^-) \ge C_{12}/(1+h)$. This yields
    \begin{equation}
    \label{pathe2}
    P^z( T_i>\min(\tau_0^-,\tau_k^+)) \ge {C_{12} \over i+1} P(R_1<-1)=:C_{13} {1\over i+1}\,.
    \end{equation}
When $i\le k/2$, (and $z \in[i,i+1)$), we notice that
	\begin{eqnarray*}
	E^z\left[\e^{U_k},\,\tau_k^+<\tau_0^-\right] &\le& E^z\left[E^{R_{\tau_{k/2}^+}}[\e^{U_k}],\, \tau_{k/2}^+ <\tau_0^-\right]\\
	&\le& C_{14} P^z(\tau_{k/2}^+ < \tau_0^-)\\
	&\le& C_{15} {i+1 \over k}
	\end{eqnarray*}
where the last two inequalities come from Lemmas \ref{overshoot} and \ref{passage}.
For $i\ge k/2$, we simply write
	$$
	E^z[\e^{U_k},\,\tau_k^+<\tau_0^-] \le \sup_{k\ge 0}E^z[e^{U_k}].
	$$
Therefore, we have for any $i\le k$,
\begin{equation}\label{pathe3}
E^z[\e^{U_k},\,\tau_k^+<\tau_0^-]  \le C_{16} {1+i \over k}.
\end{equation}

\noindent We obtain that for any integer $i$ between $1$ and $k-1$, and any $z \in I_i$,
	\begin{eqnarray}
\nonumber \qquad 	&& E^z\left[\e^{U_k}N(i),\, \tau_k^+ < \tau_0^-\right]\\
\nonumber	&\le&
	\sum_{n\ge 0} (1+n)\left(\sup_{z\in I_i} P^z(T_i < \min(\tau_k^+,\tau_0^-))\right)^n\sup_{z\in I_i}E^z\left[\e^{U_k},\,\tau_k^+< \min(\tau_0^-,T_i)\right]\\
\nonumber &=&
\left( 1 - \sup_{z\in I_i} P^z(T_i < \min(\tau_k^+,\tau_0^-))\right)^{-2}\sup_{z\in I_i}E^z\left[\e^{U_k},\,\tau_k^+< \min(\tau_0^-,T_i)\right]\\
\nonumber &\le&
 C_{13}^{-2}(i+1)^2 C_{16} {1+i\over k}\\
    &\le& C_{17} {(i+1)^3 \over k}  \label{pathgreen}
    \end{eqnarray}
by (\ref{pathe2}) and (\ref{pathe3}). We have to deal with the extreme cases $i=0$ and $i>k-1$. For $z \in I_0$,  we see that $P^z(T_i>\min(\tau_0^-, \tau_k^+)) \ge P(Y_1<-1)$, which yields by the same reasoning as before
    $$
    E^z\left[\e^{U_k}N(0),\, \tau_k^+<\tau_0^-\right] \le C_{18} {1 \over k}\,.
    $$
Similarly, ($\lfloor k\rfloor$ is the biggest integer smaller than $k$),
    $$
    E^z\left[\e^{U_k} N(\lfloor k \rfloor),\, \tau_k^+<\tau_0^-\right] \le C_{19}\,.
    $$
Therefore, (\ref{pathgreen}) still holds for any integer $i\in[0,k)$, as long as $C_{17}$ is taken large enough. By the strong Markov property, we deduce that
$$
E^0\left[\e^{U_k}N(i),\, \tau_k^+ < \tau_0^-\right] \le C_{17}P^0(T_{i}<\tau_0^-\land \tau_k^+){(i+1)^3 \over k}.
$$

\noindent This gives the following upper bound for $A$:
    \begin{equation}
  A \le C_{17}\sum_{0\le i <k} g(i+1)P^0(T_{i}<\tau_0^-\land \tau_k^+){(i+1)^3 \over k}\,.
\label{boundA}
    \end{equation}
In particular,
    $$
    A \le C_{17} {1 \over k} \sum_{0\le i <k} (i+1)^3g(i+1) = C_{20} {1 \over k}
    $$
with $C_{20}:=C_{17} \sum_{i \ge 0} (i+1)^3g(i+1)$. This proves the upper bound of (\ref{path0tok}). For the lower bound, we write (beware that $U_k\ge 0$),
    \begin{eqnarray*}
	E^0\left[\e^{U_k} \sum_{\ell=0}^{\tau_{k}^+} \e^{- R_{\ell}} (R_{\ell} + 1),\, \tau_{k}^+ <\tau_0^- \right]\nonumber
	\ge
	P^0(\tau_{k}^+ <\tau_0^- ).	
	\end{eqnarray*}
We apply (\ref{passage2}) to get the lower bound of (\ref{path0tok}).  $\Box$

\section{Some moments of the killed branching random walk}

For any $a\ge 0$ and any integer $n$, we call $Z_n(a)$ the number of particles who hit level $a$ for the first time at time $n$,
    $$
    Z_n(a) := \#\{|u|=n\,:\, \tau_0^-(u)>n-1,\, \tau_a^-(u) =n  \}
    $$
where for any $a$, $\tau_a^-(u)$ is the hitting time of $(-\infty,a)$ of the particle $u$.
We notice that particles in $Z_n(a)$ can be dead at time $n$, but their father at time $n-1$ is necessarily alive. Let also
    $$
    Z(a) := \sum_{n\ge 0} Z_n(a).
    $$
Similarly, for any $k  > a\ge 0$, and any integer $n \ge 0$, we introduce
    \begin{eqnarray*}
    Z_n(a,k) &:=& \#\{|u|=n \, : \, \tau_0^-(u) > n - 1,\, \tau_k^+(u)>n,\, \tau_a^-(u) =n   \}\,,\\
    Z(a,k)   &:=& \sum_{n\ge 0} Z_n(a,k)\,.
    \end{eqnarray*}
In words, $Z_n(a,k)$ stands for the number of particles who hit level $a$ at time $n$ and did not touch level $k$ before.

\bigskip

We denote by $S_n=X_0+X_1+ \ldots + X_n$ the random walk whose steps are distributed like $X$. We define the probability $Q^y$ as the probability which verifies for every $n$
    \begin{equation}\label{chngmeas}
    {dQ^y \over dP^y}_{| X_0,..,X_n} := {\e^{\rho (S_n - S_0)} \over \phi(\rho)^n}.
    \end{equation}

\noindent Under $Q^y$, the random walk $S_n$ is centered and starts at $y$.

\begin{proposition}
\label{moments}
We have for any $x\ge 0$, and any $a\ge 0$,
\begin{eqnarray}
\label{moment1} E^{k}[Z(a,k)] & = & \Delta_k(1) {\e^{\rho (k-a)} \over k}\,,\\
\label{moment2} E^{k}[Z(a,k)^2]  & = & \Delta_k(1){\e^{\rho (2 k- 2 a)} \over k^2}\,.
\end{eqnarray}

\noindent Besides, if $x>a \ge 0$,
\begin{equation}
\label{moment3} E^x[Z(a,k)^2]   =  \Delta_k(1)(1+x){\e^{\rho (k+x-2 a)}\over k^3}\,.
\end{equation}
\end{proposition}
{\it Proof}. Let $y$ be any real in $[0,k]$ and let $a \in [0,y]$. We observe that
    $$
    E^y[Z_n(a,k)] = b^n P^y(\tau_0^- > n - 1,\, \tau_k^+>n,\, \tau_a^- =n) \,.
    $$
The change of measure yields that
    \begin{eqnarray*}
    E^y[Z_n(a,k)] &=&   \e^{\rho y}E_Q^y[\e^{-\rho S_n},\,\tau_0^- > n - 1,\, \tau_k^+>n,\, \tau_a^- =n]\\
    &=& \e^{\rho(y-a)} E_Q^y[\e^{\rho (a-S_n)},\,\tau_0^- > n - 1,\, \tau_k^+>n,\, \tau_a^- =n]\,.
    \end{eqnarray*}
Summing over $n$ leads to
    \begin{equation}\label{firstZ}
    E^y[Z(a,k)] = \e^{\rho (y-a)}E_Q^y[\e^{\rho L_a},\,\tau_a^- < \tau_k ^+]\,.
    \end{equation}
Suppose that $y>k/2$. We observe that
	\begin{eqnarray*}
	&&E_Q^y\left[\e^{\rho L_a},\,\tau_a^- < \tau_k ^+\right] \\
    &\le& 	
	E_Q^y\left[\e^{\rho \left(a-S_{\tau_{k/2}^-}\right)},\, S_{\tau_{k/2}^-}<a\right]
	+
	E_Q^y\left[E_Q^{h}[\e^{\rho L_a}]_{h=S_{\tau_{k/2}^-}},\,\tau_{k/2}^- < \tau_k^+,\,S_{\tau_{k/2}^-}\ge a\right].
	\end{eqnarray*}
We know by Lemma \ref{overshoot} that $\sup_{\ell \le 0} E_{Q}^0[e^{\rho L_{\ell}}] \le C_{22}$. We deduce that
	$$
	E_Q^y\left[e^{\rho L_a},\,\tau_a^- < \tau_k ^+\right] \le C_{22}\left(e^{\rho(a-k/2)} + P^y\left(\tau_{k/2}^-< \tau_k^+ \right)\right)\,.	
	$$
We use Lemma \ref{passage} (applied to $R_\ell= k - S_\ell$) to see that for $k$ greater than some constant $K(a)$ (whose value may change during the proof),
	$$
	E_Q^y\left[e^{\rho L_a},\,\tau_a^- < \tau_k ^+\right]\le C_{23}{1 + k-y \over k}\,.
	$$
For $y\le k/2$, we see that
	$$
	E_Q^y\left[e^{\rho L_a },\,\tau_a^- < \tau_k ^+\right] \le
	E_Q^y\left[e^{\rho L_a}\right] \le C_{22}\,.
	$$
We deduce the existence of a constant $C_{24}$ such that for any $0\le a\le y \le k$ and any $k\ge K(a)$, we have $E_Q^y\left[e^{\rho L_a},\,\tau_a^- < \tau_k ^+\right]\le C_{24}{1 + k-y \over k}$. It yields by (\ref{firstZ}) that
	\begin{equation}\label{upunif}
	E^y[Z(a,k)] \le C_{24}e^{\rho (y-a)}{ 1 + k-y \over k}\,.
	\end{equation}
Since $L_a\ge 0$, we get  $E_Q^y[e^{\rho L_a},\,\tau_a^- < \tau_k ^+] \ge Q^y(\tau_a^- < \tau_k ^+)$. By Lemma \ref{passage},
	$$
	Q^y(\tau_a^- < \tau_k^+) \ge C_{25}{ 1+k-y\over k-a}\,.
	$$
Therefore, using (\ref{firstZ}), we get that
	\begin{equation}\label{lowunif}
	E^y[Z(a,k)] \ge {C_{25}} \, \e^{\rho (y-a)}{ 1 + k-y \over k}\,.
	\end{equation}
Equations (\ref{upunif}) and (\ref{lowunif}) give (\ref{moment1}) by taking $y=k$. We turn to the proof of  (\ref{moment2}) and (\ref{moment3}).
    \begin{eqnarray}
    \nonumber E^y[Z(a,k)^2] & = & \sum_{n\ge 0} E^y[ Z(a,k)Z_n(a,k)] \\
    \label{spine}              & = & \sum_{n\ge 0} \sum_{|u|=n} E^y[ Z(a,k),\,n= \tau_a^-(u) <\tau_k^+(u) ]\,.
    \end{eqnarray}
We decompose $Z(a,k)$ along the particle $u$ to get
    $$
    Z(a,k) = 1 + \sum_{\ell =0}^{n-1} Z^{u_{\ell}}(a,k)
    $$
where $u_{\ell}$ is the ancestor of $u$ at time $\ell$ and $Z^{u_{\ell}}(a,k)$ is the number of descendants $v$ of $u_{\ell}$ at time $n$ which are not descendants of $u_{\ell +1}$ and such that $n=\tau_a^-(v)<\tau_k^+(v)$. In particular,
    \begin{eqnarray*}
    E[Z^{u_{\ell}}(a,k)] &=& (b-1)\left(E^{S(u_{\ell})}\left[E^{S_1}[ Z(a,k)],\,S_1\in [a,k]  \right]+P^{S(u_{\ell})}(S_1 <a)\right) \\
    &=& (b-1)\left(\Delta_k(1) E^{S(u_{\ell})}\left[e^{\rho(S_1 - a)} {1+k-S_1 \over k},\, S_1\in [a,k]\right]+ P(Y_1<a-S(u_{\ell}))\right)\\
    &=& \Delta_k(1) e^{\rho (S(u_{\ell}) -a)} {1+k - S(u_{\ell}) \over k}
    \end{eqnarray*}
if $k\ge K(a)$ and $S(u_\ell)\ge a$. This decomposition leads to
    \begin{eqnarray*}
    &&E^y\left[ Z(a,k),\,n=\tau_a^-(u) <\tau_k^+(u)  \right] \\
     &=&  \Delta_k(1){e^{-\rho a} \over k}\sum_{\ell =0}^n E^y\left[e^{\rho S(u_{\ell})}(k-S(u_{\ell}) +1) ,\, n=\tau_a^-(u) <\tau_k^+(u) \right]\,.
    \end{eqnarray*}
Then equation (\ref{spine}) becomes
    \begin{eqnarray}
\nonumber    E^y[Z(a,k)^2] & = & \Delta_k(1){e^{-\rho a} \over k}\sum_{n\ge 0} b^n \sum_{\ell =0}^n E^y\left[ e^{\rho S_{\ell}} (k-S_{\ell} + 1),\, n=\tau_a^- <\tau_k^+ \right]\\
\nonumber                  & = & \Delta_k(1) {e^{\rho (y-a)}\over k}\sum_{n\ge 0} \sum_{\ell =0}^n E^y_Q\left[ e^{\rho (S_{\ell}-S_n)} (k-S_{\ell} + 1),\, n=\tau_a^- <\tau_k^+ \right] \\
\label{rappelsec}                  & = & \Delta_k(1) {e^{\rho (y-a)}\over k} E^y_Q\left[e^{-\rho S_{\tau_a^-}} \sum_{\ell=0}^{\tau_a^-} e^{\rho S_{\ell}} (k-S_{\ell} + 1) ,\, \tau_a^- <\tau_k^+ \right]
    \end{eqnarray}

\noindent where we used the change of measure from $P^y$ to $Q^y$ defined in (\ref{chngmeas}). Take $y=k$. It implies that
    \begin{eqnarray*}
   &&E^{k}[Z(a,k)^2] \\
   &=&	  \Delta_k(1) {e^{ \rho  (2k-2 a)} \over k} E^{k}_Q\left[ e^{\rho L_a}\sum_{\ell=0}^{\tau_a^-} \e^{-\rho (k-S_{\ell})} {(k-S_{\ell} + 1)},\, \tau_{a}^- <\tau_k^+ \right]\,.
    \end{eqnarray*}
We apply equation (\ref{path0tok}) of Lemma \ref{estimatepath} for the walk $R_{\ell}:= \rho(k - S_{\ell})$ to get  (\ref{moment2}). If we take $y=x$, we obtain
	\begin{eqnarray*}
   E^x[Z(a,k)^2] =  \Delta_k(1) {e^{ \rho (x+ k - 2 a)} \over k} E^{x}_Q\left[ e^{\rho L_a}\sum_{\ell=0}^{\tau_{a}^-} \e^{- \rho(k- S_{\ell})} {(k-S_{\ell} + 1)},\, \tau_{a}^- <\tau_k^+ \right]
    \end{eqnarray*}
and we apply (\ref{pathktok}) of Lemma \ref{estimatepath} to complete the proof of (\ref{moment3}). $\Box$

\section{Tail distribution of the maximum}

We are interested in large deviations of the maximum $M$ of  the branching random walk before its extinction
	$$
	M := \sup\{ S(u)\,:\, u\in \mathcal{T} \mbox{ such that } S(v)\ge 0 \; \forall \, v\le u\}\,.	
	$$
To this end, we introduce
	\begin{eqnarray*}
	H_n(k) &:=& \#\{ |u|=n\,:\, \tau_k^+(u)=n, \tau_0^-(u)>n  \}\,, \\
	H(k) &:=& \sum_{n\ge 1} H_n(k)\,.
	\end{eqnarray*}
The variable $H(k)$ is the number of particles of the branching random walk on $[0,k]$ with two killing barriers which were absorbed at level $k$.
\begin{proposition}
\label{momentsmax}
We have
\begin{eqnarray}
\label{momentsmax1} E^x[H_k] & = & \Delta_k(1) \e^{\rho (x-k)}{1+x \over k}\,,\\
\label{momentsmax2} E^x[H_k^2]  & = & \Delta_k(1)\e^{\rho (x -  k)}{1+x \over k}\,.
\end{eqnarray}
\end{proposition}

\bigskip

It shows that $H_k$ is strongly concentrated. Our result on the maximal position states as follows.
\begin{corollary}
\label{max}
The tail distribution of $M$ verifies
$$
P^x(M\ge k)=\Delta_k(1)(1+x){e^{\rho(x -k)} \over k}
$$
\end{corollary}
{\it Proof}. The corollary easily follows from the following inequalities
	$$
	P^x(M\ge k)\le E^x[H_k]
	$$
and
	$$
	P(M\ge k) = P(H_k\ge 1)\ge {E[H_k]^2\over E[H_k^2]}\,. \qquad \Box
	$$

\bigskip
	
We turn to the proof of Proposition \ref{momentsmax}. Since it is really similar to the proof of Proposition \ref{moments}, we feel free to skip some of the details. \\

\noindent {\it Proof of Proposition \ref{momentsmax}}. We verify that
	\begin{equation}\label{firstmax}
	E^x[H_k] = e^{\rho(x-k)}E^x_Q[e^{- \rho U_k},\, \tau_k^+< \tau_0^-]\,.
	\end{equation}
Since $U_k\ge 0$, we deduce that
	\begin{equation}\label{upmax}
	E^x[H_k] \le e^{\rho(x-k)}Q^x(\tau_k^+<\tau_0^-) = \Delta_k(1)e^{\rho(x-k)}{1+x\over k}\,.
	\end{equation}
On the other hand, observe that
	\begin{eqnarray*}
	E^x_Q[e^{-\rho U_k},\, \tau_k^+< \tau_0^-]
	\ge
	e^{-\rho M}Q^x(U_k<M,\, \tau_k^+< \tau_0^-)
	\end{eqnarray*}
We see that
	\begin{eqnarray*}
&&	Q^x(U_k \ge M,\, \tau_k^+< \tau_0^-)\\
	&\le&	
	Q^x\left(S_{\tau_{k/2}^+}<k,\, \tau_{k/2}^+< \tau_0^-\right)\sup_{\ell\ge 0} Q^0\left(U_{\ell} \ge M\right)
    + Q^x\left( U_{k/2}>k/2 \right)\\
	&\le& {1+x \over k}\varepsilon(M) + o(1/k)
	\end{eqnarray*}
for some $\varepsilon(M)$ which goes to zero when $M$ goes to infinity by Lemma \ref{overshoot}. Therefore
	\begin{equation}\label{lowmax}
	Q^x(U_k <M,\, \tau_k^+< \tau_0^-)\ge C_{26}{1+x \over k}
	\end{equation}
for $M$ large enough. Equations (\ref{firstmax}), (\ref{upmax}) and (\ref{lowmax}) give (\ref{momentsmax1}). We look then at the second moment of $H_k$. As before (see (\ref{rappelsec})), we can write
	\begin{eqnarray*}
	E^x[H_k^2] = \Delta_k(1) {e^{\rho(x-k)}\over k}E_Q^x\left[e^{-\rho U_k} \sum_{\ell=0}^{\tau_k^+} e^{\rho (S_{\ell} -k)}(1+S_{\ell}),\,\tau_k^+<\tau_0^-\right]\,.
	\end{eqnarray*}
We apply (\ref{pathkto0}) of Lemma \ref{estimatepath} to complete the proof. $\Box$

\section{Proof of Theorem \ref{main}}

{\it Proof of Theorem \ref{main}: lower bound}.
Let $a \in (0,x)$. We observe that
    $$
    P^x(Z>n) \ge P^x(M \ge k)P^k(Z(k,a)>n)\,.
    $$
By Proposition \ref{moments}, there exists a constant $\mu>0$ such that $E^k[Z(k,a)] \le \mu e^{\rho k}/k$ when $k$ is large enough. Let $k$ be such that $\mu e^{\rho k}/(2k)=n$. Then $k = {1 \over \rho} \ln(n) + o(\ln(n))$, and we get by Corollary \ref{max}
    $$
    P^x(M\ge k) \ge  C_{27}{(1+x)e^{\rho x} \over n\ln^2(n)}\,.
    $$
By the choice of $k$, we notice that
    $$
    P^k\left(Z(k,a)>n\right) \ge P^k\left(Z(k,a) > {E[Z(k,a)]\over 2}\right)\,.
    $$
Thus Paley-Zygmund inequality leads to
    $$
    P^k(Z(k,a) > n) \ge {1\over 4} { E^k[Z(k,a)]^2 \over E^k[Z(k,a)^2]}\,.
    $$
Proposition \ref{moments} shows then that $P^k(Z(k,a) > n) \ge C_{28}>0$. Therefore,
    $$
    P^x(Z>n) \ge C_{29}{ (1+x)e^{\rho x}\over n\ln^2(n)}
    $$
with $C_{29}=C_{27}C_{28}/4$, which proves the lower bound of the theorem. $\Box$

\bigskip
\bigskip
We turn to the proof of the upper bound. We recall that $Z(0)$ represents the number of particles who hit the domain $(-\infty,0)$.\\

\noindent {\it Proof of Theorem : upper bound}.
First, we notice that $Z(0)= 1  + (b-1) Z$. Indeed, $Z(0)$ is the number of leaves of a tree of size $Z+Z(0)$, in which any vertex has either zero or $b$ children. Therefore
    $$
    P^x(Z>n) = P^x\left( Z(0) > {n-1 \over b-1}\right)\,.
    $$
Hence it is equivalent to find an upper bound for $P^x(Z(0)>n)$. For any $k$, we have that
    \begin{eqnarray*}
    P^x(Z(0)>n) &\le& P^x(M < k,\, Z(0,k)>n) + P^x( M \ge k) \\
    &\le& P^x(Z(0,k) > n) + P^x(M \ge k)\,.
    \end{eqnarray*}
By Markov inequality, then Proposition \ref{moments}, we have
    $$
    P^x(Z(0,k) > n) \le {E^x[Z(0,k)^2] \over n^2} \le C_{30}(1+x)e^{\rho x}{ e^{\rho k} \over k^3n^2 }\,.
    $$
Therefore, by Corollary \ref{max}, we have for $k$ large enough
    $$
    P^x(Z(0)>n) \le C_{31}(1+x)e^{\rho x}\left({ e^{\rho k} \over k^3n^2 } + {e^{-\rho k} \over k}\right)\,.
    $$
Take  $k$ such that $e^{\rho k}/k=n$. We verify that
    $$
    { e^{\rho k} \over k^3n^2 } = {e^{-\rho k} \over k} = {1\over \rho^2}{1 \over n \ln^2(n)}(1+o(1))
    $$
which gives the desired upper bound. $\Box$

\bigskip
\bigskip
\bigskip

{\it Acknowledgements}. I am grateful to Zhan Shi for discussions and useful comments on the work.

\bibliographystyle{plain}
\bibliography{biblio}

\end{document}